\documentclass[11pt]{article}
\usepackage{amssymb}
\usepackage{stmaryrd}
\usepackage{indentfirst}
\usepackage{amscd}
\usepackage{amsmath}
\usepackage{latexsym}
\usepackage{amsfonts}
\usepackage{amssymb}
\usepackage{amsthm}
\usepackage{eqnarray}
\usepackage{mathrsfs}
\usepackage{ams,amsmath,amsthm,amsfonts,amssymb,bm}
\usepackage[dvips]{graphicx}
\usepackage{subfigure}

\begin{document}
\title{Characterizations of operator order \\for $k$ strictly positive operators
 \thanks
{ This work is supported by National Natural Science Fund of China
(10771011 and 11171013
).}}
\author{ Jian Shi$^{1}$ \ \ Zongsheng Gao \hskip 1cm  \\
{\small   LMIB $\&$ School of Mathematics and Systems Science,}\\
{\small Beihang University, Beijing, 100191, China}\\ }
         \date{}
         \maketitle

\maketitle \baselineskip 16pt \line(1,0){420}

\noindent{\bf Abstract.}\ \  Let $A_{i}\ (i=1, 2, \cdots , k)$ be bounded linear operators on a
Hilbert space. This paper aims to show characterizations of operator order $A_{k}\geq A_{k-1}\geq\cdots\geq A_{2}\geq A_{1}>0$ in terms of operator inequalities. Afterwards, an application of the characterizations is given to operator equalities due to Douglas's majorization and factorization theorem.\vspace{0.2cm}


 \vspace{0.2cm} \vskip   0.2cm \footnotetext[1]{Corresponding author.

E-mail
addresses: shijian@ss.buaa.edu.cn\ (J. Shi), zshgao@buaa.edu.cn\ (Z. Gao).}

\setlength{\baselineskip}{20pt}
\section{Introduction }

A capital letter (such as $T$) means a bounded linear operator on a Hilbert space $\mathcal{H}$. $T$ is said to be positive (denoted by $T\geq 0$) if $(Tx, x) \geq 0$ for all $x\in \mathcal{H} $, and $T$ is said to be strictly positive  (denoted by $T>0$) if $T$ is positive and invertible. The usual order $S\geq T$ among selfadjoint operators on $\mathcal{H}$ is defined by $(Sx, x)\geq (Tx, x)$ for all $x\in \mathcal{H}$. Let $I$ denote the indentity operator.

As an essential and historical extension of the famous L\"{o}wner-Heinz inequality: $A\geq B\geq 0 \Rightarrow A^{\alpha}\geq B^{\alpha}$ if $\alpha \in [0, 1]$, T. Furuta proved the following operator inequality in 1987.

\noindent{\bf Theorem 1.1.} (Furuta Inequality, \cite{Furuta_1987}) If $A\geq B\geq 0$, then for each $r\geq 0$,
\begin{eqnarray}
(A^{\frac r 2}A^{p}A^{\frac r 2})^{\frac 1 q}\geq (A^{\frac r 2}B^{p}A^{\frac r 2})^{\frac 1 q},\\
(B^{\frac r 2}A^{p}B^{\frac r 2})^{\frac 1 q}\geq (B^{\frac r 2}B^{p}B^{\frac r 2})^{\frac 1 q}
\end{eqnarray}
hold for $p\geq 0$ and $q\geq 1$ with $(1+r)q\geq p+r$.

\begin{figure}[h]  \centering  \includegraphics{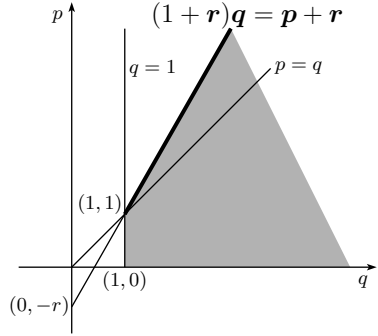} \caption{ Domain of Furuta inequality}  \label{fige1}  \end{figure}

K. Tanahashi showed that the conditions $p$ and $q$ in Figure 1 are best possible for each $r\geq 0$. See \cite{Tanahashi_1996}. It is well known that Furuta inequality has many applications. See \cite{Aluthge_1990, Aluthge_Powers, Fujii_1993, F_F, Huruya, Ito, Yuan_Gao_2008, Yuan_2010MIA, Yuan_2010OM}.

In 1995, T. Furuta showed the following theorem which interpolates Furuta inequality.

\noindent{\bf Theorem 1.2.} (Grand Furuta Inequality, \cite{Furuta_1995}) If $A\geq B\geq 0$ with $A>0$, then for each $t\in [0, 1]$ and $p\geq 1$,
\begin{eqnarray}
A^{1-t+r}\geq\{A^{\frac r 2}(A^{-{\frac t 2}}B^{p}A^{-{\frac t 2}})A^{\frac r 2}\}^{\frac {1-t+r}{(p-t)s+r}}
\end{eqnarray}
holds for $s\geq 1$ and $r\geq t$.

K. Tanahashi proved that the exponent value ${\frac {1-t+r}{(p-t)s+r}}$ of grand Furuta inequality is the best possible in \cite{Tanahashi_2000}. Afterwards, the proof was improved by T. Yamazaki and M. Fujii et al., respectively. See \cite{Yamazaki_1999} and \cite{Fujii_1999}.

In 2003, grand Furuta inequality was extended by M. Uchiyama in \cite{Uchiyama_2003} as follows:

\noindent{\bf Theorem 1.3.} (Extended Grand Furuta Inequality, \cite{Uchiyama_2003}) If $A\geq B\geq C\geq 0$ with $B>0$, then for each $t\in [0, 1]$ and $p\geq 1$,
\begin{eqnarray}
A^{1-t+r}\geq {A^{\frac r 2}(B^{-{\frac t 2}}C^{p}B^{-{\frac t 2}})^{s}A^{\frac r 2}}^{\frac {1-t+r}{(p-t)s+r}}
\end{eqnarray}
holds for $s\geq 1$ and $r\geq t$.

In 2008, grand Furuta inequality was given another extension in \cite{Furuta_2008} as follows:

\noindent{\bf Theorem 1.4.} (Extension of Furuta Inequality, \cite{Furuta_2008}) If $A\geq B\geq 0$ with $A>0$, $t\in [0, 1]$ and $p_{1}, p_{2}, \ldots, p_{2n}\geq 1$ for any natural number $n$, then the following inequality
\begin{equation}
\begin{split}
A^{1-t+r}\geq  &\ \big \{A^{\frac r 2}\big[A^{-{\frac t 2}}\cdots \big[A^{-{\frac t 2}}\{A^{\frac t 2}(A^{-{\frac t 2}}B^{p_{1}}A^{-{\frac t 2}})^{p_{2}}\\ &\ A^{\frac t 2}\}^{p_{3}}A^{-{\frac t 2}}\big]^{p_{4}}
\cdots A^{-{\frac t 2}}\big]^{p_{2n}}A^{\frac r 2}\big\}^{\frac {1-t+r}{\phi (2n)-t+r}}
\end{split}
\end{equation}
holds for $r\geq t$, where $\phi (2n)=\{\cdots[\{[(p_{1}-t)p_{2}+t]p_{3}-t\}p_{4}+t]p_{5}-\cdots -t\}p_{2n}+t$.

In 2010, C. Yang and Y. Wang showed the following theorem which interpolates extended grand Furuta inequality in \cite{Yang_2010}.

\noindent{\bf Theorem 1.5.} (Further Extension of Furuta Inequality, \cite{Yang_2010}) If $A_{2n+1}\geq A_{2n}\geq A_{2n-1}\geq \cdots\geq A_{3}\geq A_{2}\geq A_{1}\geq 0$ with $A_{2}>0$, $t_{1}, t_{2}, \ldots, t_{n-1}, t_{n} \in [0, 1]$ and $p_{1}, p_{2}, \ldots, p_{2n-1}, p_{2n}\geq 1$ for a natural number $n$, then the following inequality
\begin{equation}
\begin{split}
A^{1-t_{n}+r}_{2n+1}\geq &\  \{A^{\frac r 2}_{2n+1}[A^{-{\frac {t_{n}}{2}}}_{2n}\{A^{\frac {t_{n-1}}{2}}_{2n-1} \cdots A^{\frac {t_{2}}{2}}_{5}[A^{-{\frac {t_{2}}{2}}}_{4}\{A^{\frac {t_{1}}{2}}_{3}(A^{-{\frac {t_{1}}{2}}}_{2}A^{p_{1}}_{1}A^{-{\frac {t_{1}}{2}}}_{2})^{p_{2}}\\
&\  A^{\frac {t_{1}}{2}}_{3}\}^{p_{3}}A^{-{\frac {t_{2}}{2}}}_{4}]^{p_{4}}A^{\frac {t_{2}}{2}}_{5}
\cdots A^{\frac {t_{n-1}}{2}}_{2n-1}\}^{p_{2n-1}}A^{-{\frac {t_{n}}{2}}}_{2n}]^{p_{2n}}A^{\frac r 2}_{2n+1}\}^{\frac {1-t_{n}+r}{\psi [2n]-t_{n}+r}}
\end{split}
\end{equation}
holds for $r\geq t_{n}$, where $\psi [2n]=\{\cdots [\{[(p_{1}-t_{1})p_{2}+t_{1}]p_{3}-t_{2}\}p_{4}+t_{2}]p_{5}-\cdots -t_{n}\}p_{2n}+t_{n}$.

Recently, some beautiful results on characterizations of operator order have been shown, such as \cite{Fujii_2002}, \cite{Lin_2007} and \cite{Lin_2011}. C.-S. Lin, by using Furuta inequality, showed the characterizations of operator order for two strictly positive operators in \cite{Lin_2011}. Afterwards, he and Y. J. Cho, by using extended grand Furuta inequality, showed the characterizations of operator order for three strictly positive operators. The aim of the present paper is to show the characterizations of operator order $A_{k}\geq A_{k-1}\geq \cdots \geq A_{2}\geq A_{1}>0$ for any positive integer $k$ in terms of operator inequality via further extension of Furuta inequality. An application of the characterizations is given to operator equalities due to Douglas's majorization and factorization theorem.

\section{Main results and proofs}

\indent In this section, we show the characterizations of operator order for $k$ strictly positive operators. First, we assume that $k$ is an odd integer ($k=2n+1$).\\

\noindent{\bf Theorem 2.1.} Let $A_{1}, A_{2}, A_{3}, \cdots, A_{2n-1}, A_{2n}, A_{2n+1}$ be strictly positive operators. Then the following two assertions are equivalent.\\
(I)  $A_{2n+1}\geq A_{2n}\geq A_{2n-1}\geq \cdots \geq A_{3}\geq A_{2}\geq A_{1}$.\\
(II) If $t_{1}, t_{2}, \cdots, t_{n}\in [0, 1]$, $p_{1}, p_{2}, \cdots, p_{2n-1}, p_{2n}\geq 1$, $\psi [2n]=\{\cdots[\{[(p_{1}-t_{1})p_{2}+t_{1}]p_{3}-t_{2}\}p_{4}+t_{2}]p_{5}-\cdots -t_{n}\}p_{2n}+t_{n}$, then the following inequalities always hold for $r\geq t_{n}$:\\
\indent (II.1) $A^{r-t_{n}}_{2n+1}\geq     \Big\{A^{\frac r 2}_{2n+1}\Big[A^{-{\frac {t_{n}}{2}}}_{2n}\big\{A^{\frac {t_{n-1}}{2}}_{2n-1} \cdots A^{\frac {t_{2}}{2}}_{5}\big[A^{-{\frac {t_{2}}{2}}}_{4}\cdot\{A^{\frac {t_{1}}{2}}_{3}(A^{-{\frac {t_{1}}{2}}}_{2}A^{p_{1}}_{1}A^{-{\frac {t_{1}}{2}}}_{2})^{p_{2}} A^{\frac {t_{1}}{2}}_{3}\}^{p_{3}}\cdot \\
 A^{-{\frac {t_{2}}{2}}}_{4}\big]^{p_{4}}A^{\frac {t_{2}}{2}}_{5}
\cdots A^{\frac {t_{n-1}}{2}}_{2n-1}\big\}^{p_{2n-1}}A^{-{\frac {t_{n}}{2}}}_{2n}\Big]^{p_{2n}}A^{\frac r 2}_{2n+1}\Big\}^{\frac {r-t_{n}}{\psi [2n]-t_{n}+r}}$;\\
\indent (II.2) $A^{r-t_{n}}_{2n+1}\geq     \Big\{A^{\frac r 2}_{2n+1}\Big[A^{-{\frac {t_{n}}{2}}}_{2n+1}\big\{A^{\frac {t_{n-1}}{2}}_{2n } \cdots A^{\frac {t_{2}}{2}}_{6}\big[A^{-{\frac {t_{2}}{2}}}_{5}\cdot\{A^{\frac {t_{1}}{2}}_{4}(A^{-{\frac {t_{1}}{2}}}_{3}A^{p_{1}}_{2}A^{-{\frac {t_{1}}{2}}}_{3})^{p_{2}} A^{\frac {t_{1}}{2}}_{4}\}^{p_{3}}\cdot \\
 A^{-{\frac {t_{2}}{2}}}_{5}\big]^{p_{4}}A^{\frac {t_{2}}{2}}_{6}
\cdots A^{\frac {t_{n-1}}{2}}_{2n }\big\}^{p_{2n-1}}A^{-{\frac {t_{n}}{2}}}_{2n+1}\Big]^{p_{2n}}A^{\frac r 2}_{2n+1}\Big\}^{\frac {r-t_{n}}{\psi [2n]-t_{n}+r}}$;\\
\indent (II.3) $A^{r-t_{n}}_{2n+1}\geq     \Big\{A^{\frac r 2}_{2n+1}\Big[A^{-{\frac {t_{n}}{2}}}_{2n+1}\big\{A^{\frac {t_{n-1}}{2}}_{2n+1 } \cdots A^{\frac {t_{2}}{2}}_{7}\big[A^{-{\frac {t_{2}}{2}}}_{6}\cdot\{A^{\frac {t_{1}}{2}}_{5}(A^{-{\frac {t_{1}}{2}}}_{4}A^{p_{1}}_{3}A^{-{\frac {t_{1}}{2}}}_{4})^{p_{2}} A^{\frac {t_{1}}{2}}_{5}\}^{p_{3}} \cdot\\
 A^{-{\frac {t_{2}}{2}}}_{6}\big]^{p_{4}}A^{\frac {t_{2}}{2}}_{7}
\cdots A^{\frac {t_{n-1}}{2}}_{2n+1 }\big\}^{p_{2n-1}}A^{-{\frac {t_{n}}{2}}}_{2n+1}\Big]^{p_{2n}}A^{\frac r 2}_{2n+1}\Big\}^{\frac {r-t_{n}}{\psi [2n]-t_{n}+r}}$;\\
\indent  {  $\cdots \cdots \cdots \cdots$}\\
\indent (II.n) $A^{r-t_{n}}_{2n+1}\geq     \Big\{A^{\frac r 2}_{2n+1}\Big[A^{-{\frac {t_{n}}{2}}}_{2n+1}\big\{A^{\frac {t_{n-1}}{2}}_{2n+1 } \cdots A^{\frac {t_{2}}{2}}_{n+4}\big[A^{-{\frac {t_{2}}{2}}}_{n+3}\{A^{\frac {t_{1}}{2}}_{n+2}(A^{-{\frac {t_{1}}{2}}}_{n+1}A^{p_{1}}_{n}A^{-{\frac {t_{1}}{2}}}_{n+1})^{p_{2}} A^{\frac {t_{1}}{2}}_{n+2}\}^{p_{3}}  \\
 A^{-{\frac {t_{2}}{2}}}_{n+3}\big]^{p_{4}}A^{\frac {t_{2}}{2}}_{n+4}
\cdots A^{\frac {t_{n-1}}{2}}_{2n+1 }\big\}^{p_{2n-1}}A^{-{\frac {t_{n}}{2}}}_{2n+1}\Big]^{p_{2n}}A^{\frac r 2}_{2n+1}\Big\}^{\frac {r-t_{n}}{\psi [2n]-t_{n}+r}}$;\\
\indent (II.n+1) $A^{r-t_{n}}_{ 1}\leq     \Big\{A^{\frac r 2}_{ 1}\Big[A^{-{\frac {t_{n}}{2}}}_{1}\big\{A^{\frac {t_{n-1}}{2}}_{ 1 } \cdots A^{\frac {t_{2}}{2}}_{n-2}\big[A^{-{\frac {t_{2}}{2}}}_{n-1}\{A^{\frac {t_{1}}{2}}_{n }(A^{-{\frac {t_{1}}{2}}}_{n+1}A^{p_{1}}_{n+2}A^{-{\frac {t_{1}}{2}}}_{n+1})^{p_{2}} A^{\frac {t_{1}}{2}}_{n }\}^{p_{3}}  \\
 A^{-{\frac {t_{2}}{2}}}_{n-1}\big]^{p_{4}}A^{\frac {t_{2}}{2}}_{n-2}
\cdots A^{\frac {t_{n-1}}{2}}_{1}\big\}^{p_{2n-1}}A^{-{\frac {t_{n}}{2}}}_{1}\Big]^{p_{2n}}A^{\frac r 2}_{1}\Big\}^{\frac {r-t_{n}}{\psi [2n]-t_{n}+r}}$;\\
\indent  {  $\cdots \cdots \cdots \cdots$}\\
\indent (II.2n-2) $A^{r-t_{n}}_{ 1}\leq     \Big\{A^{\frac r 2}_{ 1}\Big[A^{-{\frac {t_{n}}{2}}}_{1}\big\{A^{\frac {t_{n-1}}{2}}_{1} \cdots A^{\frac {t_{2}}{2}}_{2n-5}\big[A^{-{\frac {t_{2}}{2}}}_{2n-4}\{A^{\frac {t_{1}}{2}}_{2n-3 }\cdot(A^{-{\frac {t_{1}}{2}}}_{2n-2}A^{p_{1}}_{2n-1}A^{-{\frac {t_{1}}{2}}}_{2n-2})^{p_{2}}\cdot \\ A^{\frac {t_{1}}{2}}_{2n-3}\}^{p_{3}}
 A^{-{\frac {t_{2}}{2}}}_{2n-4}\big]^{p_{4}}A^{\frac {t_{2}}{2}}_{2n-5}
\cdots A^{\frac {t_{n-1}}{2}}_{1}\big\}^{p_{2n-1}}A^{-{\frac {t_{n}}{2}}}_{1}\Big]^{p_{2n}}A^{\frac r 2}_{1}\Big\}^{\frac {r-t_{n}}{\psi [2n]-t_{n}+r}}$;\\
\indent (II.2n-1) $A^{r-t_{n}}_{ 1}\leq     \Big\{A^{\frac r 2}_{ 1}\Big[A^{-{\frac {t_{n}}{2}}}_{1}\big\{A^{\frac {t_{n-1}}{2}}_{2} \cdots A^{\frac {t_{2}}{2}}_{2n-4}\big[A^{-{\frac {t_{2}}{2}}}_{2n-3}\{A^{\frac {t_{1}}{2}}_{2n-2 }\cdot(A^{-{\frac {t_{1}}{2}}}_{2n-1}A^{p_{1}}_{2n}A^{-{\frac {t_{1}}{2}}}_{2n-1})^{p_{2}}\cdot \\ A^{\frac {t_{1}}{2}}_{2n-2}\}^{p_{3}}
 A^{-{\frac {t_{2}}{2}}}_{2n-3}\big]^{p_{4}}A^{\frac {t_{2}}{2}}_{2n-4}
\cdots A^{\frac {t_{n-1}}{2}}_{2}\big\}^{p_{2n-1}}A^{-{\frac {t_{n}}{2}}}_{1}\Big]^{p_{2n}}A^{\frac r 2}_{1}\Big\}^{\frac {r-t_{n}}{\psi [2n]-t_{n}+r}}$;\\
\indent (II.2n) $A^{r-t_{n}}_{ 1}\leq     \Big\{A^{\frac r 2}_{ 1}\Big[A^{-{\frac {t_{n}}{2}}}_{2}\big\{A^{\frac {t_{n-1}}{2}}_{3} \cdots A^{\frac {t_{2}}{2}}_{2n-3}\big[A^{-{\frac {t_{2}}{2}}}_{2n-2}\{A^{\frac {t_{1}}{2}}_{2n-1}\cdot(A^{-{\frac {t_{1}}{2}}}_{2n}A^{p_{1}}_{2n+1}A^{-{\frac {t_{1}}{2}}}_{2n})^{p_{2}}\cdot \\ A^{\frac {t_{1}}{2}}_{2n-1}\}^{p_{3}}
 A^{-{\frac {t_{2}}{2}}}_{2n-2}\big]^{p_{4}}A^{\frac {t_{2}}{2}}_{2n-3}
\cdots A^{\frac {t_{n-1}}{2}}_{3}\big\}^{p_{2n-1}}A^{-{\frac {t_{n}}{2}}}_{2}\Big]^{p_{2n}}A^{\frac r 2}_{1}\Big\}^{\frac {r-t_{n}}{\psi [2n]-t_{n}+r}}$.\\

\noindent {\bf Proof.} (I)$\Rightarrow $(II) Applying L\"{o}wner-Heinz inequality for ${\frac {r-t_{n}}{1-t_{n}+r}}$ to further extension of Furuta inequality, (II.1) is obtained; Replacing $A_{1}, A_{2}, A_{3}, \cdots,A_{2n-1}, A_{2n}$ by $A_{2}, A_{3}, A_{4}, \cdots, A_{2n}, A_{2n+1}$ in (II.1), respectively, (II.2) is obtained; Replacing $A_{2}, A_{3}, A_{4}, \cdots, A_{2n-1}, A_{2n}$ by $A_{3}, A_{4}, A_{5}, \cdots, A_{2n}, A_{2n+1}$ in (II.2), respectively, (II.3) is obtained. Similarly, we can obtain (II.4), (II.5), $\cdots$, (II.n).

If we replace $A_{1}, A_{2}, A_{3} \cdots, A_{2n-1}, A_{2n}, A_{2n+1}$ by $A_{2n+1}^{-1}, A_{2n}^{-1}, A_{2n-1}^{-1}, \cdots, A_{3}^{-1},$ $ A_{2}^{-1}, A_{1}^{-1}$ in (II.1), (II.2), (II.3), $\cdots$, (II.n), respectively, and take reverse, then (II.2n), (II.2n-1), (II.2n-2), $\cdots$, (II.n+1) hold.

(II)$\Rightarrow$(I) Because each $A_{i}$  is strictly positive and bounded, there exist $u_{i}$ and $v_{i}$ such that $+\infty>u_{i}I\geq A_{i}\geq v_{i}I>0$ $(i=1, 2, \cdots, 2n+1)$. If we take $p_{1}=p_{3}=p_{4}=\cdots =p_{2n}=1$, $t_{1}=t_{2}=\cdots =t_{n}=1$, $r=2$ in (II.1), then we have
\begin{equation}\tag{2.1}
\begin{split}
&\ A_{2n+1}\\ \geq &\  \big\{A_{2n+1}A_{2n}^{-{\frac 1 2}}A_{2n-1}^{\frac 1 2}\cdots A_{5}^{\frac 1 2} A_{4}^{-{\frac  1 2}}A_{3}^{\frac 1 2}(A^{-{\frac 1 2}}_{2}A_{1}A^{-{\frac 1 2}}_{2})^{p_{2}}A_{3}^{\frac 1 2}A_{4}^{-{\frac  1 2}} A_{5}^{\frac 1 2}\cdots A_{2n-1}^{\frac 1 2}A_{2n}^{-{\frac 1 2}}A_{2n+1}\big\}^{\frac 1 2}.
\end{split}\end{equation}

According to Theorem 6' in \cite{Fujii_1997}: $X\geq Y>0$ with $sI\geq X\geq tI>0 \Rightarrow {\frac {(s+t)^{2}}{4st}}X^{2}\geq Y^{2}$, we can obtain the following inequality by (2.1) and $u_{2n+1}I\geq A_{2n+1}\geq v_{2n+1}I>0$.
\begin{equation}\tag{2.2}
\begin{split}
&\ {\frac {(u_{2n+1}+v_{2n+1})^{2}}{4u_{2n+1}v_{2n+1}}}A_{2n+1}^{2}\\  \geq &\  A_{2n+1}A_{2n}^{-{\frac 1 2}}A_{2n-1}^{\frac 1 2}\cdots A_{5}^{\frac 1 2} A_{4}^{-{\frac  1 2}}A_{3}^{\frac 1 2}(A^{-{\frac 1 2}}_{2}A_{1}A^{-{\frac 1 2}}_{2})^{p_{2}}A_{3}^{\frac 1 2}A_{4}^{-{\frac  1 2}} A_{5}^{\frac 1 2}\cdots A_{2n-1}^{\frac 1 2}A_{2n}^{-{\frac 1 2}}A_{2n+1}.
\end{split}
\end{equation}
Then we have
\begin{equation}\tag{2.3}
\begin{split}
&\ {\frac {(u_{2n+1}+v_{2n+1})^{2}}{4u_{2n+1}v_{2n+1}}}\cdot {\frac {u_{2n}u_{2n-2}\cdots u_{6}u_{4}}{v_{2n-1}v_{2n-3}\cdots v_{5}v_{3}}}I\\
\geq &\ {\frac {(u_{2n+1}+v_{2n+1})^{2}}{4u_{2n+1}v_{2n+1}}}A_{3}^{-{\frac 1 2}}A_{4}^{\frac 1 2}A_{5}^{-{\frac 1 2}}\cdots A_{2n-1}^{-{\frac 1 2}}A_{2n}A_{2n-1}^{-{\frac 1 2}}\cdots A_{5}^{-{\frac 1 2}} A_{4}^{\frac 1 2} A_{3}^{-{\frac 1 2}}\\
\geq &\ (A_{2}^{-{\frac 1 2}}A_{1}A_{2}^{-{\frac 1 2}})^{p_{2}}.
\end{split}
\end{equation}
Thus,
\begin{equation}\tag{2.4}
\Big[{\frac {(u_{2n+1}+v_{2n+1})^{2}}{4u_{2n+1}v_{2n+1}}}\cdot {\frac {u_{2n}u_{2n-2}\cdots u_{6}u_{4}}{v_{2n-1}v_{2n-3}\cdots v_{5}v_{3}}}\Big]^{\frac  {1}{p_{2}}}I
\geq  A_{2}^{-{\frac 1 2}}A_{1}A_{2}^{-{\frac 1 2}}
\end{equation}
holds for any $p_{2}\geq 1$. $A_{2}\geq A_{1}$ is obtained by taking $p_{2}\rightarrow +\infty$.

Similarly, we can obtain $A_{3}\geq A_{2}, A_{4}\geq A_{3}, \cdots, A_{n+1}\geq A_{n}$ by (II.2), (II.3), $\cdots$, (II.n), respectively.

By the same setting for (2.2n), the following inequality holds according to Theorem 6 in \cite{Fujii_1997}: $X\geq Y>0 $ with $s'I\geq Y\geq t'I>0 \Rightarrow {\frac {(s'+t')^{2}}{4s't'}}X^{2}\geq Y^{2}.$
\begin{equation}\tag{2.5}
\begin{split}
&\ {\frac {(u_{1}+v_{1})^{2}}{4u_{1}v_{1}}}A_{1}A_{2}^{-{\frac 1 2}}A_{3}^{\frac 1 2}\cdots A_{2n-2}^{-{\frac 1 2}}A_{2n-1}^{\frac 1 2}(A_{2n}^{-{\frac 1 2}}A_{2n+1}A_{2n}^{-{\frac 1 2}})^{p_{2}}A_{2n-1}^{\frac 1 2}A_{2n-2}^{-{\frac 1 2}}\cdots A_{3}^{\frac 1 2}A_{2}^{-{\frac 1 2}}A_{1}\\
\geq &\ A_{1}^{2}.
\end{split}
\end{equation}
Then we have
\begin{equation}\tag{2.6}
\begin{split}
&\ A_{2n}^{-{\frac 1 2}}A_{2n+1}A_{2n}^{-{\frac 1 2}}\\ \geq &\ \Big[{\frac {4u_{1}v_{1}}{(u_{1}+v_{1})^{2}}}A_{2n-1}^{-{\frac 1 2}}A_{2n-2}^{\frac 1 2}\cdots A_{3}^{-{\frac 1 2}}A_{2}A_{3}^{-{\frac 1 2}}\cdots A_{2n-2}^{\frac 1 2}A_{2n-1}^{-{\frac 1 2}} \Big]^{\frac {1}{p_{2}}}\\ \geq &\ \Big[{\frac {4u_{1}v_{1}}{(u_{1}+v_{1})^{2}}}\cdot {\frac {v_{2}v_{4}\cdots v_{2n-2}}{u_{3}u_{5}\cdots u_{2n-1}}}\Big]^{\frac {1}{p_{2}}}I.
\end{split}
\end{equation}
$A_{2n+1}\geq A_{2n}$ is obtained by taking $p_{2}\rightarrow +\infty$ in (2.6).

Similarly, we can obtain $A_{2n}\geq A_{2n-1}, A_{2n-1}\geq A_{2n-2}, \cdots, A_{n+2}\geq A_{n+1}$ by (II.2n-1), (II.2n-2), $\cdots$, (II.n+1), respectively.                      $ \square$\\

\noindent{\bf Remark 2.1.} If $n=1$, Theorem 2.1 is the main result of \cite{Lin_2011}.\\

Next, we assume that $k$ is an even integer ($k=2n$).\\

\noindent{\bf Theorem 2.2.} Let $A_{1}, A_{2}, A_{3}, \cdots, A_{2n-1}, A_{2n}$ be strictly positive operators. Then the following two assertions are equivalent:\\
(I)  $A_{2n}\geq A_{2n-1}\geq \cdots \geq A_{3}\geq A_{2}\geq A_{1}$.\\
(II) If $t_{1}, t_{2}, \cdots, t_{n}\in [0, 1]$, $p_{1}, p_{2}, \cdots, p_{2n-1}, p_{2n}\geq 1$, $\psi [2n]=\{\cdots[\{[(p_{1}-t_{1})p_{2}+t_{1}]p_{3}-t_{2}\}p_{4}+t_{2}]p_{5}-\cdots -t_{n}\}p_{2n}+t_{n}$, then the following inequalities always hold for $r\geq t_{n}$:\\
\indent (II.1) $A^{r-t_{n}}_{2n}\geq     \Big\{A^{\frac r 2}_{2n }\Big[A^{-{\frac {t_{n}}{2}}}_{2n}\big\{A^{\frac {t_{n-1}}{2}}_{2n-1} \cdots A^{\frac {t_{2}}{2}}_{5}\big[A^{-{\frac {t_{2}}{2}}}_{4}\cdot\{A^{\frac {t_{1}}{2}}_{3}(A^{-{\frac {t_{1}}{2}}}_{2}A^{p_{1}}_{1}A^{-{\frac {t_{1}}{2}}}_{2})^{p_{2}} A^{\frac {t_{1}}{2}}_{3}\}^{p_{3}}\cdot \\
 A^{-{\frac {t_{2}}{2}}}_{4}\big]^{p_{4}}A^{\frac {t_{2}}{2}}_{5}
\cdots A^{\frac {t_{n-1}}{2}}_{2n-1}\big\}^{p_{2n-1}}A^{-{\frac {t_{n}}{2}}}_{2n}\Big]^{p_{2n}}A^{\frac r 2}_{2n }\Big\}^{\frac {r-t_{n}}{\psi [2n]-t_{n}+r}}$;\\
\indent (II.2) $A^{r-t_{n}}_{2n }\geq     \Big\{A^{\frac r 2}_{2n }\Big[A^{-{\frac {t_{n}}{2}}}_{2n }\big\{A^{\frac {t_{n-1}}{2}}_{2n } \cdots A^{\frac {t_{2}}{2}}_{6}\big[A^{-{\frac {t_{2}}{2}}}_{5}\cdot\{A^{\frac {t_{1}}{2}}_{4}(A^{-{\frac {t_{1}}{2}}}_{3}A^{p_{1}}_{2}A^{-{\frac {t_{1}}{2}}}_{3})^{p_{2}} A^{\frac {t_{1}}{2}}_{4}\}^{p_{3}}\cdot \\
 A^{-{\frac {t_{2}}{2}}}_{5}\big]^{p_{4}}A^{\frac {t_{2}}{2}}_{6}
\cdots A^{\frac {t_{n-1}}{2}}_{2n }\big\}^{p_{2n-1}}A^{-{\frac {t_{n}}{2}}}_{2n }\Big]^{p_{2n}}A^{\frac r 2}_{2n }\Big\}^{\frac {r-t_{n}}{\psi [2n]-t_{n}+r}}$;\\
\indent (II.3) $A^{r-t_{n}}_{2n }\geq     \Big\{A^{\frac r 2}_{2n }\Big[A^{-{\frac {t_{n}}{2}}}_{2n }\big\{A^{\frac {t_{n-1}}{2}}_{2n } \cdots A^{\frac {t_{2}}{2}}_{7}\big[A^{-{\frac {t_{2}}{2}}}_{6}\cdot\{A^{\frac {t_{1}}{2}}_{5}(A^{-{\frac {t_{1}}{2}}}_{4}A^{p_{1}}_{3}A^{-{\frac {t_{1}}{2}}}_{4})^{p_{2}} A^{\frac {t_{1}}{2}}_{5}\}^{p_{3}} \cdot\\
 A^{-{\frac {t_{2}}{2}}}_{6}\big]^{p_{4}}A^{\frac {t_{2}}{2}}_{7}
\cdots A^{\frac {t_{n-1}}{2}}_{2n }\big\}^{p_{2n-1}}A^{-{\frac {t_{n}}{2}}}_{2n }\Big]^{p_{2n}}A^{\frac r 2}_{2n }\Big\}^{\frac {r-t_{n}}{\psi [2n]-t_{n}+r}}$;\\
\indent  {  $\cdots \cdots \cdots \cdots$}\\
\indent (II.n) $A^{r-t_{n}}_{2n }\geq     \Big\{A^{\frac r 2}_{2n }\Big[A^{-{\frac {t_{n}}{2}}}_{2n }\big\{A^{\frac {t_{n-1}}{2}}_{2n } \cdots A^{\frac {t_{2}}{2}}_{n+4}\big[A^{-{\frac {t_{2}}{2}}}_{n+3}\{A^{\frac {t_{1}}{2}}_{n+2}(A^{-{\frac {t_{1}}{2}}}_{n+1}A^{p_{1}}_{n}A^{-{\frac {t_{1}}{2}}}_{n+1})^{p_{2}} A^{\frac {t_{1}}{2}}_{n+2}\}^{p_{3}}  \\
 A^{-{\frac {t_{2}}{2}}}_{n+3}\big]^{p_{4}}A^{\frac {t_{2}}{2}}_{n+4}
\cdots A^{\frac {t_{n-1}}{2}}_{2n }\big\}^{p_{2n-1}}A^{-{\frac {t_{n}}{2}}}_{2n }\Big]^{p_{2n}}A^{\frac r 2}_{2n }\Big\}^{\frac {r-t_{n}}{\psi [2n]-t_{n}+r}}$;\\
\indent (II.n+1) $A^{r-t_{n}}_{ 1}\leq     \Big\{A^{\frac r 2}_{ 1}\Big[A^{-{\frac {t_{n}}{2}}}_{1}\big\{A^{\frac {t_{n-1}}{2}}_{ 1 } \cdots A^{\frac {t_{2}}{2}}_{n-2}\big[A^{-{\frac {t_{2}}{2}}}_{n-1}\{A^{\frac {t_{1}}{2}}_{n }(A^{-{\frac {t_{1}}{2}}}_{n+1}A^{p_{1}}_{n+2}A^{-{\frac {t_{1}}{2}}}_{n+1})^{p_{2}} A^{\frac {t_{1}}{2}}_{n }\}^{p_{3}}  \\
 A^{-{\frac {t_{2}}{2}}}_{n-1}\big]^{p_{4}}A^{\frac {t_{2}}{2}}_{n-2}
\cdots A^{\frac {t_{n-1}}{2}}_{1}\big\}^{p_{2n-1}}A^{-{\frac {t_{n}}{2}}}_{1}\Big]^{p_{2n}}A^{\frac r 2}_{1}\Big\}^{\frac {r-t_{n}}{\psi [2n]-t_{n}+r}}$;\\
\indent  {  $\cdots \cdots \cdots \cdots$}\\
\indent (II.2n-2) $A^{r-t_{n}}_{ 1}\leq     \Big\{A^{\frac r 2}_{ 1}\Big[A^{-{\frac {t_{n}}{2}}}_{1}\big\{A^{\frac {t_{n-1}}{2}}_{1} \cdots A^{\frac {t_{2}}{2}}_{2n-5}\big[A^{-{\frac {t_{2}}{2}}}_{2n-4}\{A^{\frac {t_{1}}{2}}_{2n-3 }\cdot(A^{-{\frac {t_{1}}{2}}}_{2n-2}A^{p_{1}}_{2n-1}A^{-{\frac {t_{1}}{2}}}_{2n-2})^{p_{2}}\cdot \\ A^{\frac {t_{1}}{2}}_{2n-3}\}^{p_{3}}
 A^{-{\frac {t_{2}}{2}}}_{2n-4}\big]^{p_{4}}A^{\frac {t_{2}}{2}}_{2n-5}
\cdots A^{\frac {t_{n-1}}{2}}_{1}\big\}^{p_{2n-1}}A^{-{\frac {t_{n}}{2}}}_{1}\Big]^{p_{2n}}A^{\frac r 2}_{1}\Big\}^{\frac {r-t_{n}}{\psi [2n]-t_{n}+r}}$;\\
\indent (II.2n-1) $A^{r-t_{n}}_{ 1}\leq     \Big\{A^{\frac r 2}_{ 1}\Big[A^{-{\frac {t_{n}}{2}}}_{1}\big\{A^{\frac {t_{n-1}}{2}}_{2} \cdots A^{\frac {t_{2}}{2}}_{2n-4}\big[A^{-{\frac {t_{2}}{2}}}_{2n-3}\{A^{\frac {t_{1}}{2}}_{2n-2 }\cdot(A^{-{\frac {t_{1}}{2}}}_{2n-1}A^{p_{1}}_{2n}A^{-{\frac {t_{1}}{2}}}_{2n-1})^{p_{2}}\cdot \\ A^{\frac {t_{1}}{2}}_{2n-2}\}^{p_{3}}
 A^{-{\frac {t_{2}}{2}}}_{2n-3}\big]^{p_{4}}A^{\frac {t_{2}}{2}}_{2n-4}
\cdots A^{\frac {t_{n-1}}{2}}_{2}\big\}^{p_{2n-1}}A^{-{\frac {t_{n}}{2}}}_{1}\Big]^{p_{2n}}A^{\frac r 2}_{1}\Big\}^{\frac {r-t_{n}}{\psi [2n]-t_{n}+r}}$.\\

\noindent{\bf Proof.} Let $A_{2n+1}=A_{2n}$ in Theorem 2.1.          $\square$\\

Together with Theorem 2.1 and Theorem 2.2, we show the characterizations of operator order $A_{k}\geq A_{k-1}\geq \cdots \geq A_{2}\geq A_{1}>0$ for any positive integer $k$. For example, if $k=5$, we have the following result:

\noindent{\bf Proposition 2.1.} Let $A_{1}, A_{2}, A_{3}, A_{4}$ and $A_{5}$ be strictly positive operators. Then $A_{5}\geq A_{4}\geq A_{3}\geq A_{2}\geq A_{1}$ if and only if the following four operator inequalities
\begin{equation}\tag{2.7}
A_{5}^{r-t_{2}}\geq \Big\{A_{5}^{\frac r 2}\Big[A_{4}^{-{\frac {t_{2}} {2}}}\big(A_{3}^{ {\frac {t_{1}}{2}}}(A_{2}^{-{\frac {t_{1}}{2}}}A_{1}^{p_{1}}A_{2}^{-{\frac {t_{1}}{2}}})^{p_{2}}A_{3}^{ {\frac {t_{1}}{2}}}\big)^{p_{3}}A_{4}^{-{\frac {t_{2}} {2}}}\Big]^{p_{4}}A_{5}^{\frac r 2}\Big\}^{\frac {r-t_{2}}{\psi [4]-t_{2}+r}},
\end{equation}
\begin{equation}\tag{2.8}
A_{5}^{r-t_{2}}\geq \Big\{A_{5}^{\frac r 2}\Big[A_{5}^{-{\frac {t_{2}} {2}}}\big(A_{4}^{ {\frac {t_{1}}{2}}}(A_{3}^{-{\frac {t_{1}}{2}}}A_{2}^{p_{1}}A_{3}^{-{\frac {t_{1}}{2}}})^{p_{2}}A_{4}^{ {\frac {t_{1}}{2}}}\big)^{p_{3}}A_{5}^{-{\frac {t_{2}} {2}}}\Big]^{p_{4}}A_{5}^{\frac r 2}\Big\}^{\frac {r-t_{2}}{\psi [4]-t_{2}+r}},
\end{equation}
\begin{equation}\tag{2.9}
A_{1}^{r-t_{2}}\leq \Big\{A_{1}^{\frac r 2}\Big[A_{1}^{-{\frac {t_{2}} {2}}}\big(A_{2}^{ {\frac {t_{1}}{2}}}(A_{3}^{-{\frac {t_{1}}{2}}}A_{4}^{p_{1}}A_{3}^{-{\frac {t_{1}}{2}}})^{p_{2}}A_{2}^{ {\frac {t_{1}}{2}}}\big)^{p_{3}}A_{1}^{-{\frac {t_{2}} {2}}}\Big]^{p_{4}}A_{1}^{\frac r 2}\Big\}^{\frac {r-t_{2}}{\psi [4]-t_{2}+r}},
\end{equation}
\begin{equation}\tag{2.10}
A_{1}^{r-t_{2}}\leq \Big\{A_{1}^{\frac r 2}\Big[A_{2}^{-{\frac {t_{2}} {2}}}\big(A_{3}^{ {\frac {t_{1}}{2}}}(A_{4}^{-{\frac {t_{1}}{2}}}A_{5}^{p_{1}}A_{4}^{-{\frac {t_{1}}{2}}})^{p_{2}}A_{3}^{ {\frac {t_{1}}{2}}}\big)^{p_{3}}A_{2}^{-{\frac {t_{2}} {2}}}\Big]^{p_{4}}A_{1}^{\frac r 2}\Big\}^{\frac {r-t_{2}}{\psi [4]-t_{2}+r}}
\end{equation}
always hold for $p_{1}, p_{2}, p_{3}, p_{4}\geq 1$, $t_{1}, t_{2}\in [0, 1]$ and $r\geq t_{2}$, where $\psi [4]=\{[(p_{1}-t_{1})p_{2}+t_{1}]p_{3}-t_{2}\}p_{4}+t_{2}$.\\

\noindent{\bf Remark 2.2.} It should be mentioned that we can not obtain  $A_{5}\geq A_{4}\geq A_{3}\geq A_{2}\geq A_{1}$     only by (2.7) and (2.10).  If $A_{1}=\begin{pmatrix} 1 & 0 \\ 0 & {\frac 1 u} \end{pmatrix}$,  $A_{2}=\begin{pmatrix} 1 & 0 \\ 0 & 1 \end{pmatrix}$, $A_{3}=\begin{pmatrix} 1 & 0 \\ 0 & u \end{pmatrix}$, $A_{4}=\begin{pmatrix} u & 0 \\ 0 & 1 \end{pmatrix}$, $A_{5}=\begin{pmatrix} u+\varepsilon & 0 \\ 0 & 1 \end{pmatrix}$, where $u>1$ and $\varepsilon >0$, then the five strictly positive operators satisfy (2.7) and (2.10) without satisfying $A_{4}\geq A_{3}$.

\section{An application}

In what follows we give an application of the characterizations in Theorem 2.1 and Theorem 2.2 to operator equalities.

\noindent{\bf Theorem 3.1.} If $A_{1}, A_{2}, A_{3}, \cdots, A_{2n-1}, A_{2n}, A_{2n+1}$ are strictly positive operators, $t_{1}, t_{2}, \cdots, t_{n}\in [0, 1]$, $p_{1}, p_{2}, \cdots, p_{2n} \geq 1$, $\psi [2n]=\{\cdots [\{[(p_{1}-t_{1})p_{2}+t_{1}]p_{3}-t_{2}\}p_{4}+t_{2}]p_{5}-\cdots -t_{n}\}p_{2n}+t_{n}$, $r\geq t_{n}$, $m$ is a positive integer such that $(r-t_{n})m=\psi [2n]-t_{n}+r$ with $m\geq 2$, then the following assertions are mutually equivalent:\\
(I) $A_{2n+1}\geq A_{2n}\geq A_{2n-1}\geq \cdots \geq A_{3}\geq A_{2}\geq A_{1}$.\\
(II) The following operator inequalities hold: \\
\indent (II.1) $A^{r-t_{n}}_{2n+1}\geq     \Big\{A^{\frac r 2}_{2n+1}\Big[A^{-{\frac {t_{n}}{2}}}_{2n}\big\{A^{\frac {t_{n-1}}{2}}_{2n-1} \cdots A^{\frac {t_{2}}{2}}_{5}\big[A^{-{\frac {t_{2}}{2}}}_{4}\cdot\{A^{\frac {t_{1}}{2}}_{3}(A^{-{\frac {t_{1}}{2}}}_{2}A^{p_{1}}_{1}A^{-{\frac {t_{1}}{2}}}_{2})^{p_{2}} A^{\frac {t_{1}}{2}}_{3}\}^{p_{3}}\cdot \\
 A^{-{\frac {t_{2}}{2}}}_{4}\big]^{p_{4}}A^{\frac {t_{2}}{2}}_{5}
\cdots A^{\frac {t_{n-1}}{2}}_{2n-1}\big\}^{p_{2n-1}}A^{-{\frac {t_{n}}{2}}}_{2n}\Big]^{p_{2n}}A^{\frac r 2}_{2n+1}\Big\}^{ \frac 1 m}$;\\
\indent (II.2) $A^{r-t_{n}}_{2n+1}\geq     \Big\{A^{\frac r 2}_{2n+1}\Big[A^{-{\frac {t_{n}}{2}}}_{2n+1}\big\{A^{\frac {t_{n-1}}{2}}_{2n } \cdots A^{\frac {t_{2}}{2}}_{6}\big[A^{-{\frac {t_{2}}{2}}}_{5}\cdot\{A^{\frac {t_{1}}{2}}_{4}(A^{-{\frac {t_{1}}{2}}}_{3}A^{p_{1}}_{2}A^{-{\frac {t_{1}}{2}}}_{3})^{p_{2}} A^{\frac {t_{1}}{2}}_{4}\}^{p_{3}}\cdot \\
 A^{-{\frac {t_{2}}{2}}}_{5}\big]^{p_{4}}A^{\frac {t_{2}}{2}}_{6}
\cdots A^{\frac {t_{n-1}}{2}}_{2n }\big\}^{p_{2n-1}}A^{-{\frac {t_{n}}{2}}}_{2n+1}\Big]^{p_{2n}}A^{\frac r 2}_{2n+1}\Big\}^{ \frac 1 m}$;\\
\indent (II.3) $A^{r-t_{n}}_{2n+1}\geq     \Big\{A^{\frac r 2}_{2n+1}\Big[A^{-{\frac {t_{n}}{2}}}_{2n+1}\big\{A^{\frac {t_{n-1}}{2}}_{2n+1 } \cdots A^{\frac {t_{2}}{2}}_{7}\big[A^{-{\frac {t_{2}}{2}}}_{6}\cdot\{A^{\frac {t_{1}}{2}}_{5}(A^{-{\frac {t_{1}}{2}}}_{4}A^{p_{1}}_{3}A^{-{\frac {t_{1}}{2}}}_{4})^{p_{2}} A^{\frac {t_{1}}{2}}_{5}\}^{p_{3}} \cdot\\
 A^{-{\frac {t_{2}}{2}}}_{6}\big]^{p_{4}}A^{\frac {t_{2}}{2}}_{7}
\cdots A^{\frac {t_{n-1}}{2}}_{2n+1 }\big\}^{p_{2n-1}}A^{-{\frac {t_{n}}{2}}}_{2n+1}\Big]^{p_{2n}}A^{\frac r 2}_{2n+1}\Big\}^{ \frac 1 m}$;\\
\indent  {  $\cdots \cdots \cdots \cdots$}\\
\indent (II.n) $A^{r-t_{n}}_{2n+1}\geq     \Big\{A^{\frac r 2}_{2n+1}\Big[A^{-{\frac {t_{n}}{2}}}_{2n+1}\big\{A^{\frac {t_{n-1}}{2}}_{2n+1 } \cdots A^{\frac {t_{2}}{2}}_{n+4}\big[A^{-{\frac {t_{2}}{2}}}_{n+3}\{A^{\frac {t_{1}}{2}}_{n+2}(A^{-{\frac {t_{1}}{2}}}_{n+1}A^{p_{1}}_{n}A^{-{\frac {t_{1}}{2}}}_{n+1})^{p_{2}} A^{\frac {t_{1}}{2}}_{n+2}\}^{p_{3}}  \\
 A^{-{\frac {t_{2}}{2}}}_{n+3}\big]^{p_{4}}A^{\frac {t_{2}}{2}}_{n+4}
\cdots A^{\frac {t_{n-1}}{2}}_{2n+1 }\big\}^{p_{2n-1}}A^{-{\frac {t_{n}}{2}}}_{2n+1}\Big]^{p_{2n}}A^{\frac r 2}_{2n+1}\Big\}^{ \frac 1 m}$;\\
\indent (II.n+1) $A^{r-t_{n}}_{ 1}\leq     \Big\{A^{\frac r 2}_{ 1}\Big[A^{-{\frac {t_{n}}{2}}}_{1}\big\{A^{\frac {t_{n-1}}{2}}_{ 1 } \cdots A^{\frac {t_{2}}{2}}_{n-2}\big[A^{-{\frac {t_{2}}{2}}}_{n-1}\{A^{\frac {t_{1}}{2}}_{n }(A^{-{\frac {t_{1}}{2}}}_{n+1}A^{p_{1}}_{n+2}A^{-{\frac {t_{1}}{2}}}_{n+1})^{p_{2}} A^{\frac {t_{1}}{2}}_{n }\}^{p_{3}}  \\
 A^{-{\frac {t_{2}}{2}}}_{n-1}\big]^{p_{4}}A^{\frac {t_{2}}{2}}_{n-2}
\cdots A^{\frac {t_{n-1}}{2}}_{1}\big\}^{p_{2n-1}}A^{-{\frac {t_{n}}{2}}}_{1}\Big]^{p_{2n}}A^{\frac r 2}_{1}\Big\}^{ \frac 1 m}$;\\
\indent  {  $\cdots \cdots \cdots \cdots$}\\
\indent (II.2n-2) $A^{r-t_{n}}_{ 1}\leq     \Big\{A^{\frac r 2}_{ 1}\Big[A^{-{\frac {t_{n}}{2}}}_{1}\big\{A^{\frac {t_{n-1}}{2}}_{1} \cdots A^{\frac {t_{2}}{2}}_{2n-5}\big[A^{-{\frac {t_{2}}{2}}}_{2n-4}\{A^{\frac {t_{1}}{2}}_{2n-3 }\cdot(A^{-{\frac {t_{1}}{2}}}_{2n-2}A^{p_{1}}_{2n-1}A^{-{\frac {t_{1}}{2}}}_{2n-2})^{p_{2}}\cdot \\ A^{\frac {t_{1}}{2}}_{2n-3}\}^{p_{3}}
 A^{-{\frac {t_{2}}{2}}}_{2n-4}\big]^{p_{4}}A^{\frac {t_{2}}{2}}_{2n-5}
\cdots A^{\frac {t_{n-1}}{2}}_{1}\big\}^{p_{2n-1}}A^{-{\frac {t_{n}}{2}}}_{1}\Big]^{p_{2n}}A^{\frac r 2}_{1}\Big\}^{ \frac 1 m}$;\\
\indent (II.2n-1) $A^{r-t_{n}}_{ 1}\leq     \Big\{A^{\frac r 2}_{ 1}\Big[A^{-{\frac {t_{n}}{2}}}_{1}\big\{A^{\frac {t_{n-1}}{2}}_{2} \cdots A^{\frac {t_{2}}{2}}_{2n-4}\big[A^{-{\frac {t_{2}}{2}}}_{2n-3}\{A^{\frac {t_{1}}{2}}_{2n-2 }\cdot(A^{-{\frac {t_{1}}{2}}}_{2n-1}A^{p_{1}}_{2n}A^{-{\frac {t_{1}}{2}}}_{2n-1})^{p_{2}}\cdot \\ A^{\frac {t_{1}}{2}}_{2n-2}\}^{p_{3}}
 A^{-{\frac {t_{2}}{2}}}_{2n-3}\big]^{p_{4}}A^{\frac {t_{2}}{2}}_{2n-4}
\cdots A^{\frac {t_{n-1}}{2}}_{2}\big\}^{p_{2n-1}}A^{-{\frac {t_{n}}{2}}}_{1}\Big]^{p_{2n}}A^{\frac r 2}_{1}\Big\}^{ \frac 1 m}$;\\
\indent (II.2n) $A^{r-t_{n}}_{ 1}\leq     \Big\{A^{\frac r 2}_{ 1}\Big[A^{-{\frac {t_{n}}{2}}}_{2}\big\{A^{\frac {t_{n-1}}{2}}_{3} \cdots A^{\frac {t_{2}}{2}}_{2n-3}\big[A^{-{\frac {t_{2}}{2}}}_{2n-2}\{A^{\frac {t_{1}}{2}}_{2n-1}\cdot(A^{-{\frac {t_{1}}{2}}}_{2n}A^{p_{1}}_{2n+1}A^{-{\frac {t_{1}}{2}}}_{2n})^{p_{2}}\cdot \\ A^{\frac {t_{1}}{2}}_{2n-1}\}^{p_{3}}
 A^{-{\frac {t_{2}}{2}}}_{2n-2}\big]^{p_{4}}A^{\frac {t_{2}}{2}}_{2n-3}
\cdots A^{\frac {t_{n-1}}{2}}_{3}\big\}^{p_{2n-1}}A^{-{\frac {t_{n}}{2}}}_{2}\Big]^{p_{2n}}A^{\frac r 2}_{1}\Big\}^{ \frac 1 m}$.\\
(III) There exists strictly positive operators $S_{1}, S_{2}, S_{3}, \cdots, S_{2n-2}, S_{2n-1}, S_{2n }$ satisfying the following operator equalities, respectively, where each $S_{i}$ $(i=1, 2, \cdots, 2n)$ is unique with $\|S_{i}\|\leq 1$.\\
\indent(III.1) $A_{2n+1}^{-{\frac {t_{n}}{2}}}S_{1}(A_{2n+1}^{r-t_{n}}S_{1})^{m-1}A_{2n+1}^{-{\frac {t_{n}}{2}}}=A_{2n+1}^{-{\frac {t_{n}}{2}}}(S_{1}A_{2n+1}^{r-t_{n}})^{m-1}S_{1}A_{2n+1}^{-{\frac {t_{n}}{2}}}=\\ \Big[A^{-{\frac {t_{n}}{2}}}_{2n}\big\{A^{\frac {t_{n-1}}{2}}_{2n-1} \cdots A^{\frac {t_{2}}{2}}_{5}\big[A^{-{\frac {t_{2}}{2}}}_{4} \{A^{\frac {t_{1}}{2}}_{3}(A^{-{\frac {t_{1}}{2}}}_{2}A^{p_{1}}_{1}A^{-{\frac {t_{1}}{2}}}_{2})^{p_{2}} A^{\frac {t_{1}}{2}}_{3}\}^{p_{3}}
 A^{-{\frac {t_{2}}{2}}}_{4}\big]^{p_{4}}A^{\frac {t_{2}}{2}}_{5}
\cdots A^{\frac {t_{n-1}}{2}}_{2n-1}\big\}^{p_{2n-1}}A^{-{\frac {t_{n}}{2}}}_{2n}\Big]^{p_{2n}}$;\\
\indent(III.2) $A_{2n+1}^{-{\frac {t_{n}}{2}}}S_{2}(A_{2n+1}^{r-t_{n}}S_{2})^{m-1}A_{2n+1}^{-{\frac {t_{n}}{2}}}=A_{2n+1}^{-{\frac {t_{n}}{2}}}(S_{2}A_{2n+1}^{r-t_{n}})^{m-1}S_{2}A_{2n+1}^{-{\frac {t_{n}}{2}}}=\\  \Big[A^{-{\frac {t_{n}}{2}}}_{2n+1}\big\{A^{\frac {t_{n-1}}{2}}_{2n } \cdots A^{\frac {t_{2}}{2}}_{6}\big[A^{-{\frac {t_{2}}{2}}}_{5} \{A^{\frac {t_{1}}{2}}_{4}(A^{-{\frac {t_{1}}{2}}}_{3}A^{p_{1}}_{2}A^{-{\frac {t_{1}}{2}}}_{3})^{p_{2}} A^{\frac {t_{1}}{2}}_{4}\}^{p_{3}}
 A^{-{\frac {t_{2}}{2}}}_{5}\big]^{p_{4}}A^{\frac {t_{2}}{2}}_{6}
\cdots A^{\frac {t_{n-1}}{2}}_{2n }\big\}^{p_{2n-1}}A^{-{\frac {t_{n}}{2}}}_{2n+1}\Big]^{p_{2n}}$;\\
\indent(III.3) $A_{2n+1}^{-{\frac {t_{n}}{2}}}S_{3}(A_{2n+1}^{r-t_{n}}S_{3})^{m-1}A_{2n+1}^{-{\frac {t_{n}}{2}}}=A_{2n+1}^{-{\frac {t_{n}}{2}}}(S_{3}A_{2n+1}^{r-t_{n}})^{m-1}S_{3}A_{2n+1}^{-{\frac {t_{n}}{2}}}=\\ \Big[A^{-{\frac {t_{n}}{2}}}_{2n+1}\big\{A^{\frac {t_{n-1}}{2}}_{2n+1 } \cdots A^{\frac {t_{2}}{2}}_{7}\big[A^{-{\frac {t_{2}}{2}}}_{6} \{A^{\frac {t_{1}}{2}}_{5}(A^{-{\frac {t_{1}}{2}}}_{4}A^{p_{1}}_{3}A^{-{\frac {t_{1}}{2}}}_{4})^{p_{2}} A^{\frac {t_{1}}{2}}_{5}\}^{p_{3}}
 A^{-{\frac {t_{2}}{2}}}_{6}\big]^{p_{4}}A^{\frac {t_{2}}{2}}_{7}
\cdots A^{\frac {t_{n-1}}{2}}_{2n+1 }\big\}^{p_{2n-1}}A^{-{\frac {t_{n}}{2}}}_{2n+1}\Big]^{p_{2n}}$;\\
\indent  {  $\cdots \cdots \cdots \cdots$}\\
\indent(III.n) $A_{2n+1}^{-{\frac {t_{n}}{2}}}S_{n}(A_{2n+1}^{r-t_{n}}S_{n})^{m-1}A_{2n+1}^{-{\frac {t_{n}}{2}}}=A_{2n+1}^{-{\frac {t_{n}}{2}}}(S_{n}A_{2n+1}^{r-t_{n}})^{m-1}S_{n}A_{2n+1}^{-{\frac {t_{n}}{2}}}=\\ \Big[A^{-{\frac {t_{n}}{2}}}_{2n+1}\big\{A^{\frac {t_{n-1}}{2}}_{2n+1 } \cdots A^{\frac {t_{2}}{2}}_{n+4}\big[A^{-{\frac {t_{2}}{2}}}_{n+3}\{A^{\frac {t_{1}}{2}}_{n+2}(A^{-{\frac {t_{1}}{2}}}_{n+1}A^{p_{1}}_{n}A^{-{\frac {t_{1}}{2}}}_{n+1})^{p_{2}} A^{\frac {t_{1}}{2}}_{n+2}\}^{p_{3}}
 A^{-{\frac {t_{2}}{2}}}_{n+3}\big]^{p_{4}}A^{\frac {t_{2}}{2}}_{n+4}
\cdots A^{\frac {t_{n-1}}{2}}_{2n+1 }\big\}^{p_{2n-1}}A^{-{\frac {t_{n}}{2}}}_{2n+1}\Big]^{p_{2n}}$;\\
\indent(III.n+1) $A_{1}^{-{\frac {t_{n}}{2}}}S_{n+1}^{-1}(A_{1}^{r-t_{n}}S_{n+1}^{-1})^{m-1}A_{1}^{-{\frac {t_{n}}{2}}}=A_{1}^{-{\frac {t_{n}}{2}}}(S_{n+1}^{-1}A_{1}^{r-t_{n}})^{m-1}S_{n+1}^{-1}A_{1}^{-{\frac {t_{n}}{2}}}=\\\Big[A^{-{\frac {t_{n}}{2}}}_{1}\big\{A^{\frac {t_{n-1}}{2}}_{ 1 } \cdots A^{\frac {t_{2}}{2}}_{n-2}\big[A^{-{\frac {t_{2}}{2}}}_{n-1}\{A^{\frac {t_{1}}{2}}_{n }(A^{-{\frac {t_{1}}{2}}}_{n+1}A^{p_{1}}_{n+2}A^{-{\frac {t_{1}}{2}}}_{n+1})^{p_{2}} A^{\frac {t_{1}}{2}}_{n }\}^{p_{3}}
 A^{-{\frac {t_{2}}{2}}}_{n-1}\big]^{p_{4}}A^{\frac {t_{2}}{2}}_{n-2}
\cdots A^{\frac {t_{n-1}}{2}}_{1}\big\}^{p_{2n-1}}A^{-{\frac {t_{n}}{2}}}_{1}\Big]^{p_{2n}}$;\\
\indent  {  $\cdots \cdots \cdots \cdots$}\\
\indent(III.2n-2) $ A_{1}^{-{\frac {t_{n}}{2}}}S_{2n-2}^{-1}(A_{1}^{r-t_{n}}S_{2n-2}^{-1})^{m-1}A_{1}^{-{\frac {t_{n}}{2}}}=A_{1}^{-{\frac {t_{n}}{2}}}(S_{2n-2}^{-1}A_{1}^{r-t_{n}})^{m-1}S_{2n-2}^{-1}A_{1}^{-{\frac {t_{n}}{2}}}=\\ \Big[A^{-{\frac {t_{n}}{2}}}_{1} \big\{A^{\frac {t_{n-1}}{2}}_{1} \cdots \big[A^{-{\frac {t_{2}}{2}}}_{2n-4}\{A^{\frac {t_{1}}{2}}_{2n-3 } (A^{-{\frac {t_{1}}{2}}}_{2n-2}A^{p_{1}}_{2n-1}A^{-{\frac {t_{1}}{2}}}_{2n-2})^{p_{2}}  A^{\frac {t_{1}}{2}}_{2n-3}\}^{p_{3}}
 A^{-{\frac {t_{2}}{2}}}_{2n-4}\big]^{p_{4}}
\cdots A^{\frac {t_{n-1}}{2}}_{1}\big\}^{p_{2n-1}}  A^{-{\frac {t_{n}}{2}}}_{1}\Big]^{p_{2n}}$;\\
\indent(III.2n-1) $ A_{1}^{-{\frac {t_{n}}{2}}}S_{2n-1}^{-1}(A_{1}^{r-t_{n}}S_{2n-1}^{-1})^{m-1}A_{1}^{-{\frac {t_{n}}{2}}}=A_{1}^{-{\frac {t_{n}}{2}}}(S_{2n-1}^{-1}A_{1}^{r-t_{n}})^{m-1}S_{2n-1}^{-1}A_{1}^{-{\frac {t_{n}}{2}}}=\\ \Big[A^{-{\frac {t_{n}}{2}}}_{1}\big\{A^{\frac {t_{n-1}}{2}}_{2} \cdots  \big[A^{-{\frac {t_{2}}{2}}}_{2n-3}\{A^{\frac {t_{1}}{2}}_{2n-2 } (A^{-{\frac {t_{1}}{2}}}_{2n-1}A^{p_{1}}_{2n}A^{-{\frac {t_{1}}{2}}}_{2n-1})^{p_{2}} A^{\frac {t_{1}}{2}}_{2n-2}\}^{p_{3}}
 A^{-{\frac {t_{2}}{2}}}_{2n-3}\big]^{p_{4}}
\cdots A^{\frac {t_{n-1}}{2}}_{2}\big\}^{p_{2n-1}}A^{-{\frac {t_{n}}{2}}}_{1}\Big]^{p_{2n}}$;\\
\indent(III.2n) $A_{1}^{-{\frac {t_{n}}{2}}}S_{2n }^{-1}(A_{1}^{r-t_{n}}S_{2n }^{-1})^{m-1}A_{1}^{-{\frac {t_{n}}{2}}}=A_{1}^{-{\frac {t_{n}}{2}}}(S_{2n }^{-1}A_{1}^{r-t_{n}})^{m-1}S_{2n }^{-1}A_{1}^{-{\frac {t_{n}}{2}}}=\\ \Big[A^{-{\frac {t_{n}}{2}}}_{2}\big\{A^{\frac {t_{n-1}}{2}}_{3} \cdots  \big[A^{-{\frac {t_{2}}{2}}}_{2n-2}\{A^{\frac {t_{1}}{2}}_{2n-1} (A^{-{\frac {t_{1}}{2}}}_{2n}A^{p_{1}}_{2n+1}A^{-{\frac {t_{1}}{2}}}_{2n})^{p_{2}} A^{\frac {t_{1}}{2}}_{2n-1}\}^{p_{3}}
 A^{-{\frac {t_{2}}{2}}}_{2n-2}\big]^{p_{4}}
\cdots A^{\frac {t_{n-1}}{2}}_{3}\big\}^{p_{2n-1}}A^{-{\frac {t_{n}}{2}}}_{2}\Big]^{p_{2n}}$.\\

\noindent{\bf Proof.} Because (II)$\Leftrightarrow$(III) holds obviously by Theorem 2.1, we only need to prove that (II)$\Leftrightarrow$(III).

Firstly, let us prove that (II.1)$\Rightarrow$(III.1). We recall Douglas's majorization and factorization theorem in \cite{Douglas_1966}: $SS^{\ast}\leq \lambda^{2}TT^{\ast}\Leftrightarrow $ there exists an operator $Q$ s.t. $TQ=S$, where $\|Q\|^{2}=inf\{\mu: SS^{\ast}\leq \mu TT^{\ast}\}$.\\
By (II.1), there exists an operator $E_{1}$ with $\|E_{1}\|\leq 1$ such that
\begin{equation}
\tag{3.1}
\begin{split}
&\ A_{2n+1}^{\frac {r-t_{n}}{2}}E_{1}=E_{1}^{\ast}A_{2n+1}^{\frac {r-t_{n}}{2}}\\
=&\ \Big\{A^{\frac r 2}_{2n+1}\Big[A^{-{\frac {t_{n}}{2}}}_{2n}\big\{A^{\frac {t_{n-1}}{2}}_{2n-1} \cdots  \{A^{\frac {t_{1}}{2}}_{3}(A^{-{\frac {t_{1}}{2}}}_{2}A^{p_{1}}_{1}A^{-{\frac {t_{1}}{2}}}_{2})^{p_{2}} A^{\frac {t_{1}}{2}}_{3}\}^{p_{3}}
\cdots A^{\frac {t_{n-1}}{2}}_{2n-1}\big\}^{p_{2n-1}}A^{-{\frac {t_{n}}{2}}}_{2n}\Big]^{p_{2n}}A^{\frac r 2}_{2n+1}\Big\}^{ \frac {1} {2m}}.
\end{split}
\end{equation}
Taking $S_{1}=E_{1}E_{1}^{\ast}$, we have
\begin{equation}
\tag{3.2}
\begin{split}
&\ A_{2n+1}^{\frac {r-t_{n}}{2}}S_{1}A_{2n+1}^{\frac {r-t_{n}}{2}}\\
=&\ \Big\{A^{\frac r 2}_{2n+1}\Big[A^{-{\frac {t_{n}}{2}}}_{2n}\big\{A^{\frac {t_{n-1}}{2}}_{2n-1} \cdots  \{A^{\frac {t_{1}}{2}}_{3}(A^{-{\frac {t_{1}}{2}}}_{2}A^{p_{1}}_{1}A^{-{\frac {t_{1}}{2}}}_{2})^{p_{2}} A^{\frac {t_{1}}{2}}_{3}\}^{p_{3}}
\cdots A^{\frac {t_{n-1}}{2}}_{2n-1}\big\}^{p_{2n-1}}A^{-{\frac {t_{n}}{2}}}_{2n}\Big]^{p_{2n}}A^{\frac r 2}_{2n+1}\Big\}^{ \frac {1} { m}}.
\end{split}
\end{equation}
According to (3.2) and $S_{1}=E_{1}E_{1}^{\ast}$, $S_{1}$ is unique and strictly positive with $\|S_{1}\|\leq 1$.  (3.2) also implies that
\begin{equation}
\tag{3.3}
\begin{split}
&\ (A_{2n+1}^{\frac {r-t_{n}}{2}}S_{1}A_{2n+1}^{\frac {r-t_{n}}{2}})^{m}
=  A_{2n+1}^{ {\frac {r-t_{n}}{2}}}S_{1}(A_{2n+1}^{r-t_{n}}S_{1})^{m-1}A_{2n+1}^{ {\frac {r-t_{n}}{2}}}
=  A_{2n+1}^{ {\frac {r-t_{n}}{2}}}(S_{1}A_{2n+1}^{r-t_{n}})^{m-1}S_{1}A_{2n+1}^{ {\frac {r-t_{n}}{2}}}\\
=&\ A^{\frac r 2}_{2n+1}\Big[A^{-{\frac {t_{n}}{2}}}_{2n}\big\{A^{\frac {t_{n-1}}{2}}_{2n-1} \cdots  \{A^{\frac {t_{1}}{2}}_{3}(A^{-{\frac {t_{1}}{2}}}_{2}A^{p_{1}}_{1}A^{-{\frac {t_{1}}{2}}}_{2})^{p_{2}} A^{\frac {t_{1}}{2}}_{3}\}^{p_{3}}
\cdots A^{\frac {t_{n-1}}{2}}_{2n-1}\big\}^{p_{2n-1}}A^{-{\frac {t_{n}}{2}}}_{2n}\Big]^{p_{2n}}A^{\frac r 2}_{2n+1}.
\end{split}
\end{equation}
Then (III.1) holds by (3.3).

Secondly, we prove that (III.1)$\Rightarrow$ (II.1). By (III.1),
\begin{equation*}
\begin{split}
&\ \Big\{A^{\frac r 2}_{2n+1}\Big[A^{-{\frac {t_{n}}{2}}}_{2n}\big\{A^{\frac {t_{n-1}}{2}}_{2n-1} \cdots  \{A^{\frac {t_{1}}{2}}_{3}(A^{-{\frac {t_{1}}{2}}}_{2}A^{p_{1}}_{1}A^{-{\frac {t_{1}}{2}}}_{2})^{p_{2}} A^{\frac {t_{1}}{2}}_{3}\}^{p_{3}}
\cdots A^{\frac {t_{n-1}}{2}}_{2n-1}\big\}^{p_{2n-1}}A^{-{\frac {t_{n}}{2}}}_{2n}\Big]^{p_{2n}}A^{\frac r 2}_{2n+1}\Big\}^{ \frac 1 m}\\
= &\ \big\{A_{2n+1}^{ {\frac {r-t_{n}}{2}}}S_{1}(A_{2n+1}^{r-t_{n}}S_{1})^{m-1}A_{2n+1}^{ {\frac {r-t_{n}}{2}}}\big\}^{\frac 1 m}\\
=&\ \big\{(A_{2n+1}^{ {\frac {r-t_{n}}{2}}}S_{1}A_{2n+1}^{ {\frac {r-t_{n}}{2}}})\cdot(A_{2n+1}^{ {\frac {r-t_{n}}{2}}}S_{1}A_{2n+1}^{ {\frac {r-t_{n}}{2}}})\cdots (A_{2n+1}^{ {\frac {r-t_{n}}{2}}}S_{1}A_{2n+1}^{ {\frac {r-t_{n}}{2}}})\big\}^{\frac 1 m}\\
=&\ A_{2n+1}^{ {\frac {r-t_{n}}{2}}}S_{1}A_{2n+1}^{ {\frac {r-t_{n}}{2}}}\\
\leq &\ A_{2n+1}^{r-t_{n}}.
\end{split}
\end{equation*}
The inequality follows form the fact that $S_{1}\leq \|S_{1}\|I\leq I$, and then (II.1) holds.

By using the same method above, we can prove that (II.2)$\Leftrightarrow$ (III.2), (II.3)$\Leftrightarrow$ (III.3), $\cdots$, (II.n)$\Leftrightarrow$ (III.n), respectively.

Next, we show that (II.2n)$\Leftrightarrow$(III.2n). Notice that for two strictly positive operators $S$ and $T$, $S\geq T$ if and only if $T^{-1}\geq S^{-1}$. Then (II.2n) is equivalent to
\begin{equation}
\tag{3.4}
\begin{split}
A^{-(r-t_{n})}_{ 1}\geq  &\   \Big\{A^{-{\frac r 2}}_{ 1}\Big[A^{-{\frac {t_{n}}{2}}}_{2}\big\{A^{\frac {t_{n-1}}{2}}_{3} \cdots A^{\frac {t_{2}}{2}}_{2n-3}\big[A^{-{\frac {t_{2}}{2}}}_{2n-2}\{A^{\frac {t_{1}}{2}}_{2n-1} (A^{-{\frac {t_{1}}{2}}}_{2n}A^{p_{1}}_{2n+1}A^{-{\frac {t_{1}}{2}}}_{2n})^{p_{2}}  \\&\ A^{\frac {t_{1}}{2}}_{2n-1}\}^{p_{3}}
 A^{-{\frac {t_{2}}{2}}}_{2n-2}\big]^{p_{4}}A^{\frac {t_{2}}{2}}_{2n-3}
\cdots A^{\frac {t_{n-1}}{2}}_{3}\big\}^{p_{2n-1}}A^{-{\frac {t_{n}}{2}}}_{2}\Big]^{-{p_{2n}}}A^{-{\frac r 2}}_{1}\Big\}^{ \frac 1 m}.
\end{split}
\end{equation}
The proof that (3.4)$\Leftrightarrow$(III.2n) is similar to the proof of that (II.1)$\Leftrightarrow$(III.1), so we omit it here.

Repeat the method above, we can prove that (II.2n-1)$\Leftrightarrow$(III.2n-1), (II.2n-2)$\Leftrightarrow$(III.2n-2), $\cdots $, (II.n+1)$\Leftrightarrow$(III.n+1).                  $\square$\\

\noindent{\bf Theorem 3.2.} If $A_{1}, A_{2}, A_{3}, \cdots, A_{2n-1}, A_{2n}$ are strictly positive operators, $t_{1}, t_{2}, \cdots$, $t_{n}\in [0, 1]$, $p_{1}, p_{2}, \cdots, p_{2n} \geq 1$, $\psi [2n]=\{\cdots [\{[(p_{1}-t_{1})p_{2}+t_{1}]p_{3}-t_{2}\}p_{4}+t_{2}]p_{5}-\cdots -t_{n}\}p_{2n}+t_{n}$, $r\geq t_{n}$, $m$ is a positive integer such that $(r-t_{n})m=\psi [2n]-t_{n}+r$ with $m\geq 2$, then the following assertions are mutually equivalent:\\
(I) $ A_{2n}\geq A_{2n-1}\geq \cdots \geq A_{3}\geq A_{2}\geq A_{1}$.\\
(II) The following operator inequalities hold: \\
\indent (II.1) $A^{r-t_{n}}_{2n}\geq     \Big\{A^{\frac r 2}_{2n }\Big[A^{-{\frac {t_{n}}{2}}}_{2n}\big\{A^{\frac {t_{n-1}}{2}}_{2n-1} \cdots A^{\frac {t_{2}}{2}}_{5}\big[A^{-{\frac {t_{2}}{2}}}_{4}\cdot\{A^{\frac {t_{1}}{2}}_{3}(A^{-{\frac {t_{1}}{2}}}_{2}A^{p_{1}}_{1}A^{-{\frac {t_{1}}{2}}}_{2})^{p_{2}} A^{\frac {t_{1}}{2}}_{3}\}^{p_{3}}\cdot \\
 A^{-{\frac {t_{2}}{2}}}_{4}\big]^{p_{4}}A^{\frac {t_{2}}{2}}_{5}
\cdots A^{\frac {t_{n-1}}{2}}_{2n-1}\big\}^{p_{2n-1}}A^{-{\frac {t_{n}}{2}}}_{2n}\Big]^{p_{2n}}A^{\frac r 2}_{2n }\Big\}^{\frac 1 m }$;\\
\indent (II.2) $A^{r-t_{n}}_{2n }\geq     \Big\{A^{\frac r 2}_{2n }\Big[A^{-{\frac {t_{n}}{2}}}_{2n }\big\{A^{\frac {t_{n-1}}{2}}_{2n } \cdots A^{\frac {t_{2}}{2}}_{6}\big[A^{-{\frac {t_{2}}{2}}}_{5}\cdot\{A^{\frac {t_{1}}{2}}_{4}(A^{-{\frac {t_{1}}{2}}}_{3}A^{p_{1}}_{2}A^{-{\frac {t_{1}}{2}}}_{3})^{p_{2}} A^{\frac {t_{1}}{2}}_{4}\}^{p_{3}}\cdot \\
 A^{-{\frac {t_{2}}{2}}}_{5}\big]^{p_{4}}A^{\frac {t_{2}}{2}}_{6}
\cdots A^{\frac {t_{n-1}}{2}}_{2n }\big\}^{p_{2n-1}}A^{-{\frac {t_{n}}{2}}}_{2n }\Big]^{p_{2n}}A^{\frac r 2}_{2n }\Big\}^{\frac 1 m }$;\\
\indent (II.3) $A^{r-t_{n}}_{2n }\geq     \Big\{A^{\frac r 2}_{2n }\Big[A^{-{\frac {t_{n}}{2}}}_{2n }\big\{A^{\frac {t_{n-1}}{2}}_{2n } \cdots A^{\frac {t_{2}}{2}}_{7}\big[A^{-{\frac {t_{2}}{2}}}_{6}\cdot\{A^{\frac {t_{1}}{2}}_{5}(A^{-{\frac {t_{1}}{2}}}_{4}A^{p_{1}}_{3}A^{-{\frac {t_{1}}{2}}}_{4})^{p_{2}} A^{\frac {t_{1}}{2}}_{5}\}^{p_{3}} \cdot\\
 A^{-{\frac {t_{2}}{2}}}_{6}\big]^{p_{4}}A^{\frac {t_{2}}{2}}_{7}
\cdots A^{\frac {t_{n-1}}{2}}_{2n }\big\}^{p_{2n-1}}A^{-{\frac {t_{n}}{2}}}_{2n }\Big]^{p_{2n}}A^{\frac r 2}_{2n }\Big\}^{\frac 1 m }$;\\
\indent  {  $\cdots \cdots \cdots \cdots$}\\
\indent (II.n) $A^{r-t_{n}}_{2n }\geq     \Big\{A^{\frac r 2}_{2n }\Big[A^{-{\frac {t_{n}}{2}}}_{2n }\big\{A^{\frac {t_{n-1}}{2}}_{2n } \cdots A^{\frac {t_{2}}{2}}_{n+4}\big[A^{-{\frac {t_{2}}{2}}}_{n+3}\{A^{\frac {t_{1}}{2}}_{n+2}(A^{-{\frac {t_{1}}{2}}}_{n+1}A^{p_{1}}_{n}A^{-{\frac {t_{1}}{2}}}_{n+1})^{p_{2}} A^{\frac {t_{1}}{2}}_{n+2}\}^{p_{3}}  \\
 A^{-{\frac {t_{2}}{2}}}_{n+3}\big]^{p_{4}}A^{\frac {t_{2}}{2}}_{n+4}
\cdots A^{\frac {t_{n-1}}{2}}_{2n }\big\}^{p_{2n-1}}A^{-{\frac {t_{n}}{2}}}_{2n }\Big]^{p_{2n}}A^{\frac r 2}_{2n }\Big\}^{\frac 1 m }$;\\
\indent (II.n+1) $A^{r-t_{n}}_{ 1}\leq     \Big\{A^{\frac r 2}_{ 1}\Big[A^{-{\frac {t_{n}}{2}}}_{1}\big\{A^{\frac {t_{n-1}}{2}}_{ 1 } \cdots A^{\frac {t_{2}}{2}}_{n-2}\big[A^{-{\frac {t_{2}}{2}}}_{n-1}\{A^{\frac {t_{1}}{2}}_{n }(A^{-{\frac {t_{1}}{2}}}_{n+1}A^{p_{1}}_{n+2}A^{-{\frac {t_{1}}{2}}}_{n+1})^{p_{2}} A^{\frac {t_{1}}{2}}_{n }\}^{p_{3}}  \\
 A^{-{\frac {t_{2}}{2}}}_{n-1}\big]^{p_{4}}A^{\frac {t_{2}}{2}}_{n-2}
\cdots A^{\frac {t_{n-1}}{2}}_{1}\big\}^{p_{2n-1}}A^{-{\frac {t_{n}}{2}}}_{1}\Big]^{p_{2n}}A^{\frac r 2}_{1}\Big\}^{\frac 1 m }$;\\
\indent  {  $\cdots \cdots \cdots \cdots$}\\
\indent (II.2n-2) $A^{r-t_{n}}_{ 1}\leq     \Big\{A^{\frac r 2}_{ 1}\Big[A^{-{\frac {t_{n}}{2}}}_{1}\big\{A^{\frac {t_{n-1}}{2}}_{1} \cdots A^{\frac {t_{2}}{2}}_{2n-5}\big[A^{-{\frac {t_{2}}{2}}}_{2n-4}\{A^{\frac {t_{1}}{2}}_{2n-3 }\cdot(A^{-{\frac {t_{1}}{2}}}_{2n-2}A^{p_{1}}_{2n-1}A^{-{\frac {t_{1}}{2}}}_{2n-2})^{p_{2}}\cdot \\ A^{\frac {t_{1}}{2}}_{2n-3}\}^{p_{3}}
 A^{-{\frac {t_{2}}{2}}}_{2n-4}\big]^{p_{4}}A^{\frac {t_{2}}{2}}_{2n-5}
\cdots A^{\frac {t_{n-1}}{2}}_{1}\big\}^{p_{2n-1}}A^{-{\frac {t_{n}}{2}}}_{1}\Big]^{p_{2n}}A^{\frac r 2}_{1}\Big\}^{\frac 1 m }$;\\
\indent (II.2n-1) $A^{r-t_{n}}_{ 1}\leq     \Big\{A^{\frac r 2}_{ 1}\Big[A^{-{\frac {t_{n}}{2}}}_{1}\big\{A^{\frac {t_{n-1}}{2}}_{2} \cdots A^{\frac {t_{2}}{2}}_{2n-4}\big[A^{-{\frac {t_{2}}{2}}}_{2n-3}\{A^{\frac {t_{1}}{2}}_{2n-2 }\cdot(A^{-{\frac {t_{1}}{2}}}_{2n-1}A^{p_{1}}_{2n}A^{-{\frac {t_{1}}{2}}}_{2n-1})^{p_{2}}\cdot \\ A^{\frac {t_{1}}{2}}_{2n-2}\}^{p_{3}}
 A^{-{\frac {t_{2}}{2}}}_{2n-3}\big]^{p_{4}}A^{\frac {t_{2}}{2}}_{2n-4}
\cdots A^{\frac {t_{n-1}}{2}}_{2}\big\}^{p_{2n-1}}A^{-{\frac {t_{n}}{2}}}_{1}\Big]^{p_{2n}}A^{\frac r 2}_{1}\Big\}^{\frac 1 m }$.\\
(III) There exists strictly positive operators $S_{1}, S_{2}, S_{3}, \cdots, S_{2n-2}, S_{2n-1}$ satisfying the following operator equalities, respectively, where each $S_{i}$ $(i=1, 2, \cdots, 2n-1)$ is unique with $\|S_{i}\|\leq 1$.\\
\indent(III.1) $A_{2n}^{-{\frac {t_{n}}{2}}}S_{1}(A_{2n}^{r-t_{n}}S_{1})^{m-1}A_{2n}^{-{\frac {t_{n}}{2}}}=A_{2n}^{-{\frac {t_{n}}{2}}}(S_{1}A_{2n}^{r-t_{n}})^{m-1}S_{1}A_{2n}^{-{\frac {t_{n}}{2}}}=\\ \Big[A^{-{\frac {t_{n}}{2}}}_{2n}\big\{A^{\frac {t_{n-1}}{2}}_{2n-1} \cdots A^{\frac {t_{2}}{2}}_{5}\big[A^{-{\frac {t_{2}}{2}}}_{4} \{A^{\frac {t_{1}}{2}}_{3}(A^{-{\frac {t_{1}}{2}}}_{2}A^{p_{1}}_{1}A^{-{\frac {t_{1}}{2}}}_{2})^{p_{2}} A^{\frac {t_{1}}{2}}_{3}\}^{p_{3}}
 A^{-{\frac {t_{2}}{2}}}_{4}\big]^{p_{4}}A^{\frac {t_{2}}{2}}_{5}
\cdots A^{\frac {t_{n-1}}{2}}_{2n-1}\big\}^{p_{2n-1}}A^{-{\frac {t_{n}}{2}}}_{2n}\Big]^{p_{2n}}$;\\
\indent(III.2) $A_{2n}^{-{\frac {t_{n}}{2}}}S_{2}(A_{2n}^{r-t_{n}}S_{2})^{m-1}A_{2n}^{-{\frac {t_{n}}{2}}}=A_{2n}^{-{\frac {t_{n}}{2}}}(S_{2}A_{2n}^{r-t_{n}})^{m-1}S_{2}A_{2n}^{-{\frac {t_{n}}{2}}}=\\  \Big[A^{-{\frac {t_{n}}{2}}}_{2n}\big\{A^{\frac {t_{n-1}}{2}}_{2n } \cdots A^{\frac {t_{2}}{2}}_{6}\big[A^{-{\frac {t_{2}}{2}}}_{5} \{A^{\frac {t_{1}}{2}}_{4}(A^{-{\frac {t_{1}}{2}}}_{3}A^{p_{1}}_{2}A^{-{\frac {t_{1}}{2}}}_{3})^{p_{2}} A^{\frac {t_{1}}{2}}_{4}\}^{p_{3}}
 A^{-{\frac {t_{2}}{2}}}_{5}\big]^{p_{4}}A^{\frac {t_{2}}{2}}_{6}
\cdots A^{\frac {t_{n-1}}{2}}_{2n }\big\}^{p_{2n-1}}A^{-{\frac {t_{n}}{2}}}_{2n}\Big]^{p_{2n}}$;\\
\indent(III.3) $A_{2n}^{-{\frac {t_{n}}{2}}}S_{3}(A_{2n}^{r-t_{n}}S_{3})^{m-1}A_{2n}^{-{\frac {t_{n}}{2}}}=A_{2n}^{-{\frac {t_{n}}{2}}}(S_{3}A_{2n}^{r-t_{n}})^{m-1}S_{3}A_{2n}^{-{\frac {t_{n}}{2}}}=\\ \Big[A^{-{\frac {t_{n}}{2}}}_{2n }\big\{A^{\frac {t_{n-1}}{2}}_{2n } \cdots A^{\frac {t_{2}}{2}}_{7}\big[A^{-{\frac {t_{2}}{2}}}_{6} \{A^{\frac {t_{1}}{2}}_{5}(A^{-{\frac {t_{1}}{2}}}_{4}A^{p_{1}}_{3}A^{-{\frac {t_{1}}{2}}}_{4})^{p_{2}} A^{\frac {t_{1}}{2}}_{5}\}^{p_{3}}
 A^{-{\frac {t_{2}}{2}}}_{6}\big]^{p_{4}}A^{\frac {t_{2}}{2}}_{7}
\cdots A^{\frac {t_{n-1}}{2}}_{2n }\big\}^{p_{2n-1}}A^{-{\frac {t_{n}}{2}}}_{2n }\Big]^{p_{2n}}$;\\
\indent  {  $\cdots \cdots \cdots \cdots$}\\
\indent(III.n) $A_{2n}^{-{\frac {t_{n}}{2}}}S_{n}(A_{2n}^{r-t_{n}}S_{n})^{m-1}A_{2n}^{-{\frac {t_{n}}{2}}}=A_{2n}^{-{\frac {t_{n}}{2}}}(S_{n}A_{2n}^{r-t_{n}})^{m-1}S_{n}A_{2n}^{-{\frac {t_{n}}{2}}}=\\ \Big[A^{-{\frac {t_{n}}{2}}}_{2n }\big\{A^{\frac {t_{n-1}}{2}}_{2n } \cdots A^{\frac {t_{2}}{2}}_{n+4}\big[A^{-{\frac {t_{2}}{2}}}_{n+3}\{A^{\frac {t_{1}}{2}}_{n+2}(A^{-{\frac {t_{1}}{2}}}_{n+1}A^{p_{1}}_{n}A^{-{\frac {t_{1}}{2}}}_{n+1})^{p_{2}} A^{\frac {t_{1}}{2}}_{n+2}\}^{p_{3}}
 A^{-{\frac {t_{2}}{2}}}_{n+3}\big]^{p_{4}}A^{\frac {t_{2}}{2}}_{n+4}
\cdots A^{\frac {t_{n-1}}{2}}_{2n  }\big\}^{p_{2n-1}}A^{-{\frac {t_{n}}{2}}}_{2n }\Big]^{p_{2n}}$;\\
\indent(III.n+1) $A_{1}^{-{\frac {t_{n}}{2}}}S_{n+1}^{-1}(A_{1}^{r-t_{n}}S_{n+1}^{-1})^{m-1}A_{1}^{-{\frac {t_{n}}{2}}}=A_{1}^{-{\frac {t_{n}}{2}}}(S_{n+1}^{-1}A_{1}^{r-t_{n}})^{m-1}S_{n+1}^{-1}A_{1}^{-{\frac {t_{n}}{2}}}=\\\Big[A^{-{\frac {t_{n}}{2}}}_{1}\big\{A^{\frac {t_{n-1}}{2}}_{ 1 } \cdots A^{\frac {t_{2}}{2}}_{n-2}\big[A^{-{\frac {t_{2}}{2}}}_{n-1}\{A^{\frac {t_{1}}{2}}_{n }(A^{-{\frac {t_{1}}{2}}}_{n+1}A^{p_{1}}_{n+2}A^{-{\frac {t_{1}}{2}}}_{n+1})^{p_{2}} A^{\frac {t_{1}}{2}}_{n }\}^{p_{3}}
 A^{-{\frac {t_{2}}{2}}}_{n-1}\big]^{p_{4}}A^{\frac {t_{2}}{2}}_{n-2}
\cdots A^{\frac {t_{n-1}}{2}}_{1}\big\}^{p_{2n-1}}A^{-{\frac {t_{n}}{2}}}_{1}\Big]^{p_{2n}}$;\\
\indent  {  $\cdots \cdots \cdots \cdots$}\\
\indent(III.2n-2) $ A_{1}^{-{\frac {t_{n}}{2}}}S_{2n-2}^{-1}(A_{1}^{r-t_{n}}S_{2n-2}^{-1})^{m-1}A_{1}^{-{\frac {t_{n}}{2}}}=A_{1}^{-{\frac {t_{n}}{2}}}(S_{2n-2}^{-1}A_{1}^{r-t_{n}})^{m-1}S_{2n-2}^{-1}A_{1}^{-{\frac {t_{n}}{2}}}=\\ \Big[A^{-{\frac {t_{n}}{2}}}_{1} \big\{A^{\frac {t_{n-1}}{2}}_{1} \cdots \big[A^{-{\frac {t_{2}}{2}}}_{2n-4}\{A^{\frac {t_{1}}{2}}_{2n-3 } (A^{-{\frac {t_{1}}{2}}}_{2n-2}A^{p_{1}}_{2n-1}A^{-{\frac {t_{1}}{2}}}_{2n-2})^{p_{2}}  A^{\frac {t_{1}}{2}}_{2n-3}\}^{p_{3}}
 A^{-{\frac {t_{2}}{2}}}_{2n-4}\big]^{p_{4}}
\cdots A^{\frac {t_{n-1}}{2}}_{1}\big\}^{p_{2n-1}}  A^{-{\frac {t_{n}}{2}}}_{1}\Big]^{p_{2n}}$;\\
\indent(III.2n-1) $ A_{1}^{-{\frac {t_{n}}{2}}}S_{2n-1}^{-1}(A_{1}^{r-t_{n}}S_{2n-1}^{-1})^{m-1}A_{1}^{-{\frac {t_{n}}{2}}}=A_{1}^{-{\frac {t_{n}}{2}}}(S_{2n-1}^{-1}A_{1}^{r-t_{n}})^{m-1}S_{2n-1}^{-1}A_{1}^{-{\frac {t_{n}}{2}}}=\\ \Big[A^{-{\frac {t_{n}}{2}}}_{1}\big\{A^{\frac {t_{n-1}}{2}}_{2} \cdots  \big[A^{-{\frac {t_{2}}{2}}}_{2n-3}\{A^{\frac {t_{1}}{2}}_{2n-2 } (A^{-{\frac {t_{1}}{2}}}_{2n-1}A^{p_{1}}_{2n}A^{-{\frac {t_{1}}{2}}}_{2n-1})^{p_{2}} A^{\frac {t_{1}}{2}}_{2n-2}\}^{p_{3}}
 A^{-{\frac {t_{2}}{2}}}_{2n-3}\big]^{p_{4}}
\cdots A^{\frac {t_{n-1}}{2}}_{2}\big\}^{p_{2n-1}}A^{-{\frac {t_{n}}{2}}}_{1}\Big]^{p_{2n}}$;\\

\noindent{\bf Proof.} Let $A_{2n+1}=A_{2n}$ in Theorem 3.1.             $\square$\\

Together with Theorem 3.1 and Theorem 3.2, we give an application of the characterizations of $A_{k}\geq A_{k-1}\geq \cdots \geq A_{2}\geq A_{1}>0$ to operator equalities for any positive integer $k$.

\begin{center}

\end{center}
\end{document}